\documentclass[10pt]{article}
\usepackage{amsmath}

\usepackage[english]{babel}
\usepackage{amsmath}
\usepackage{amsthm,amssymb,latexsym}

\newtheorem{prop}{Proposition}

\newtheorem{theorem}[prop]{Theorem}
\theoremstyle{definition}

\newtheorem{ex}{Example}
\newtheorem{remark}[prop]{Remark}

\title {Crucial Words and the Complexity of Some Extremal Problems for Sets of Prohibited Words}
\author {A. Evdokimov and S. Kitaev}
\begin{document}
\maketitle
\begin{abstract} 
We introduced the notation of a set of prohibitions and give definitions of a complete set and a crucial word with respect to a given set of prohibitions. We consider 3 particular sets which appear in different areas of mathematics and for each of them examine the length of a crucial word. One of these sets is proved to be incomplete. The problem of determining lengths of words that are free from a set of prohibitions is shown to be NP-complete, although the related problem of whether or not a given set of prohibitions is complete is known to be effectively solvable. 
\end{abstract}

\section{Introduction and Background}

In defining or characterising sets of objects in discrete mathematics, "languages of prohibitions" are often used to define a class of objects by listing those prohibited subobjects that are not contained in the objects of the class. To this end the notion of a subobject is defined in different ways. The notion depends on the set under consideration. These sets are subwords for partially bounded languages, subgraphs for families of graphs and so on. One of the classes of interest that have appeared and are considered in different areas of mathematics is the class of nonrecurrent symbolic sequences defined by prohibiting strong periodicity in them, or, to be more exact, by prohibiting the repetition of subwords in these symbolic sequences, for example of type $XX$.

In this paper we consider 3 types of "prohibitions" connected with a generalisation of the notion of nonrecurrent symbolic sequences, and for each of these sets we consider the structure of crucial words and find their lengths. In Section 5 we investigate the problem of determining lengths of words that are free from any given set of prohibitions. We show that this problem is NP-complete although the related problem whether or not a given set of prohibitions is complete is known to be effectively solvable. 

Let ${\bf A}=\{ a_1, \ldots ,a_n\}$ be an alphabet of $n$ letters. A {\em word} in the alphabet ${\bf A}$ is a finite sequence of letters of the alphabet. Any $i$ consecutive letters of a word $X$ generate a {\em subword} of length $i$.  If $X$ is a subword of a word $Y$, we write $X \subseteq Y$. 

The set ${\bf A}^*$ is the set of all the words in the alphabet ${\bf A}$. Let ${\bf S} \subseteq {\bf A}^*$. Then ${\bf S}$ is called a {\em set of prohibited words} or a {\em set of prohibitions}. A word that does not contain any words from ${\bf S}$ as its subwords is said to be {\em free from} ${\bf S}$. The set of all words that are free from ${\bf S}$ is denoted by $\widehat{{\bf S}}$.

\begin{ex} Let ${\bf A} = \{ a,\ b \}$. The set of prohibitions is ${\bf S} = \{ aa,\ ba \}$. The word $abbb$ is in $\widehat{{\bf S}}$. \end{ex}

If there exists a $k \in {\bf N}$ such that the length of any word in $\widehat{{\bf S}}$ is less than $k$, then ${\bf S}$ is called a {\em complete} set. 

\begin{ex} ${\bf A} = \{ 1,\ 2,\ 3,\ 4 \}$. The set of prohibitions is $${\bf S} = \{ 123,\ 13,\ 14,\ 11,\ 22,\ 33,\ 44\ \}.$$ 
Then ${\bf S}$ is incomplete, since the word $\underbrace{124124 \ldots 124}_{3k}$ is in $\widehat{{\bf S}}$ for any $k$. \end{ex}

\begin{ex} ${\bf A} = \{ 1, 2, 3 \}$. The set of prohibitions is $${\bf S} = \{ 12,\ 23,\ 31,\ 32,\ 11,\ 22,\ 33 \}.$$ 
It is easy to check that ${\bf S}$ is complete. \end{ex}

A word $X \in \widehat{{\bf S}}$ is called a {\em crucial} word (with respect to ${\bf S}$), if the word $Xa_i$ contains a prohibited subword for any letter $a_i \in {\bf A}$. This means that $Xa_i$ has the structure $BB_ia_i$, where $B$ is some word and $B_ia_i \in {\bf S}$. The subword $B_i$ is called the {\em $i$-ending} of crucial word $X$. If for each letter of the alphabet we consider minimal $i$-ending (with respect to inclusion) we obtain a system of included $i$-endings, which we will use to investigate crucial words.

\begin{ex} ${\bf A} = \{ a,\ b,\ c \}$. The set of prohibitions is ${\bf S} = \{ aa,\ cab,\ acac\}$. The word $abaca$ is crucial with respect to ${\bf S}$.\end{ex} 

A crucial word of minimal (maximal) length, if it exists, is called a {\em minimal} ({\em maximal}) crucial word.

\begin{ex} ${\bf A} = \{ a,\ b,\ c \}$. The set of prohibitions is ${\bf S} = \{ aa,\ cab,\ acac\}$. The word $aca$ is a minimal crucial word with respect to ${\bf S}$. There do not exist any maximal crucial words, since the word $\underbrace{b \ldots b}_{k}aca$ is crucial  for all $k \in {\bf N}$.\end{ex}

Let $L_{min}({\bf S})\ \ (L_{max}({\bf S}))$ denote of the length of a minimal\ (maximal) crucial word with respect to ${\bf S}$.

In this paper we consider three sets of prohibitions denoted ${\bf S}_1^n$, ${\bf S}_2^n$, ${\bf S}_3^{n, k}$. Here we use $n$ for indicating the number of letters of the alphabet under consideration and $k$ is a natural number. 

We now give the definitions of these sets:

${\bf S}_1^n=\{ XX\ |\ X \in {\bf A}^* \}$, that is, we prohibit the repetition of two equal consecutive subwords.

${\bf S}_2^n=\{ XY\ |\ \overline{\nu}(X)=\overline{\nu}(Y) \}$, where $\overline{\nu}($X$)=(\nu_1(X), \ldots ,\nu_n(X))$ is the {\em content} vector of $X$, in which $\nu_i(X)$ is the number of occurrences of the letter $a_i$ in $X$. That is, we prohibit the repetition of two consecutive subwords of the same content.

${\bf S}_3^{n, k}=\{ XY\ |\ d(X,Y) \le k, |X|=|Y| \ge k+1, k \in {\bf N} \}$, where $d(X,Y)$ is the number of letters in which the words $X$ and $Y$ differ (Hamming metric) and $|X|$ is the length of the word $X$. That is we prohibit any two consecutive subwords of the length greater then $k$ such that the number of positions in which these words differ is less then or equal to $k$.

The proofs of the theorems in this paper consist of the constructions of extremal crucial words and of the proofs of their optimality, i. e. the lower bound for $L_{min}({\bf S})$ and the upper bound for $L_{max}({\bf S})$.

\section{The Set of Prohibitions ${\bf S}_1^n$}

\begin{theorem}\label{theorem1} We have $$L_{min}({\bf S}_1^n)=2^n-1.$$ \end{theorem} 

\begin{proof} We define a crucial word $X$ by induction:

$$X_1=a_1,\  X_i=X_{i-1}a_iX_{i-1},\  X=X_n.$$

From this construction it follows that $|X|=2^n-1$. We will prove that $X$ is a minimal crucial word with respect to ${\bf S}_1^n$.

Let $U$ be an arbitrary minimal crucial word. We show that $U$ coincides with the word $X$ up to a permutation of letters in {$\bf A$}.

From the definition of a crucial word it follows that in the word $Ua_i$ there is a prohibited word of the form $B_ia_iB_ia_i$, where $B_i$ is a certain word and $B_ia_iB_ia_i$ is the ending of the word $Ua_i$ (the ending may coincide with $Ua_i$). In this case the $i$-ending is the subword $B_ia_iB_i$. Let $\ell_i=B_ia_iB_i$.

We assume that $\ell_1 \subset \ell_2 \subset \ldots \subset \ell_n$, since we can make such ordering by permuting the letters of the alphabet, which obviously does not affect the cruciality and minimality of a word.

Note that the minimal crucial word $U$ has the form
$$U=B_na_nB_n=B_na_nY_na_1,$$
where $Y_n$ is a certain word. Actually, if on the right of $B_na_nB_n$ there is a certain word, then it contradicts the minimality of a crucial word, and if instead of $a_1$ there stands $a_k$ ($k>1$) then it contradicts $\ell_1 \subset \ell_k$.

We show that $\ell_{n-1}$ coincides with $B_n$. We have $\ell_{n-1}=B_{n-1}a_{n-1}B_{n-1}$ and let $a_nB_n$ be a subword of $\ell_{n-1}$. Now $\ell_{n-1}$ has the form $Ka_nPa_{n-1}Ka_nP$, where $Ka_nP=B_{n-1}$), but then 
$$\ell_n=Pa_{n-1}Ka_nPa_{n-1}Ka_nP,\ \mbox{where}\ Pa_{n-1}Ka_nP=B_n,$$
and the word $U$ contains the prohibited subword $a_nPa_nP$. This can not be the case. It means that $\ell_{n-1}$ is a subword of the word $B_n$, and the word $U$ has the form:
$$U=\ell_n=Z_n\ell_{n-1}a_nZ_n\ell_{n-1},$$
where $Z_n$ is a certain word. Since we explore a minimal crucial word, we have $Z_n= \emptyset$, and then $B_n=\ell_{n-1}$. In the same way we can show that $B_i=\ell_{i-1}$ for each $i=2, \ldots ,n-1$ and $B_1= \emptyset$. 

Hence the structure of a minimal crucial word $U$ coincides with that of the word $X$ as required. \end{proof}

\begin{remark} From the proof of Theorem~\ref{theorem1} it follows that the word $X$ is the unique minimal crucial word to within a transposition of the letters of the alphabet ${\bf A}$. \end{remark}

\section{The Set of Prohibitions ${\bf S}_2^n$}

\begin{prop}\label{prop} A minimal crucial (with respect to ${\bf S}_2^n$) word can not have three letters, each of which appears twice in the word.  \end{prop}

\begin{proof} Since the proposition is obviously true for $|{\bf A}| = 1, 2, 3$, we will consider the case $|{\bf A}| \ge 4$. 

Let $X$ be a minimal crucial word, and suppose the system of included $i$-endings for it is $\ell_1 \subset \ell_2 \subset \ldots \subset \ell_n=X$. Suppose the letters $a_{i_1}$, $a_{i_2}$, $a_{i_3}$ occur twice in $X$ and that $i_1 < i_2 < i_3 < n$ (the fact that $i_1 , i_2 , i_3$ do not equal $n$ follows from the fact that $a_n$ must occur an odd number of times). 

When we pass from $\ell_{i_3-1}$ to $\ell_{i_3}$ ($\ell_{i_3-1}$ is determined, since there are $i_1$, $i_2 < i_3$) there must appear a letter $a_{i_3}$, and when we pass from $\ell_{i_3}$ to $\ell_{i_3+1}$ ($\ell_{i_3+1}$ is determined, since $i_3 < n$) there must appear one more letter $a_{i_3}$; Hence, since there are two letters $a_{i_3}$ in $X$, there are no letters $a_{i_3}$ for $2 < j < i_3$ in $\ell_j$ whence there are no letters $a_{i_3}$ in the $X$ to the left of $\ell_{i_2}$ (both letters $a_{i_3}$ lie to the left respecting of $\ell_{i_2}$).

Obviously, the letter $a_{i_1}$ must be in $\ell_{i_1}$. The second letter $a_{i_1}$ appears when we pass from $\ell_{i_1}$ to $\ell_{i_2}$. Since there are only two letters $a_{i_1}$, there are no letters $a_{i_1}$ in the word $X$ to the left of $\ell_{i_2}$.

If we write the letter $a_{i_3+1}$ to the right of the word $X$ we obtain a prohibited word (a word from ${\bf S}_2^n$). Words from ${\bf S}_2^n$ are divided into two parts which have the same content. Obviously, the letters $a_{i_3}$ must be in different parts of the prohibited word, and letters $a_{i_1}$ must be in different parts of the same word which is impossible, since the letters $a_{i_3}$ lie strictly to the left of $a_{i_1}$, and this contradicts the assumption. \end{proof} 

\begin{remark} From the proof of proposition 1 we have that if letters $a_i$ and $a_j$ occur twice in a word $X$ (in which $\ell_1 \subset \ell_2 \subset \ldots \subset \ell_n=X$), then either $i = j + 1$ or $j = i+1$. \end{remark}

\begin{theorem} For any $n>2$ we have $$L_{min}({\bf S}_2^n)=4n-7.$$ \end{theorem}

\begin{proof} Note that a natural approach to the construction of a crucial word is possible. It consists of an algorithm of step-by-step optimisation: We ascribe to a crucial word of an $n$-letter alphabet a minimum number of letters to obtain a crucial word of an ($n+1$)-letter alphabet. 

The algorithm can be written recursively in the following way:

$X_n=B_{n-1}a_nB_{n-1}$ 

$B_{n-1}=B_{n-3}a_{n-1}B_{n-3}$ 

$B_1=a_1,B_2=a_2,B_{-1}=B_0=X_0= \emptyset .$

\

Some initial values when implementing the algorithm are:

\

$X_1=a_1,$ 

$X_2=a_1a_2a_1,$ 

$X_3=a_2a_3a_1a_2a_1,$ 

$X_4=a_1a_3a_1a_4a_2a_3a_1a_2a_1,$ 

$X_5=a_2a_4a_2a_5a_1a_3a_1a_4a_2a_3a_1a_2a_1.$
 
\

This is an algorithm by which the minimal crucial word $X_n$ for the set of prohibitions ${\bf S}_1^n$ can be built. For ${\bf S}_2^n$ such a construction gives an upper bound of the form $exp(n/2)$, or, to be more exact, $$(3-(n \bmod 2))2^{ \lfloor \frac{n+1}{2} \rfloor}-3.$$

We now give an upper bound that is a linear function. 

We introduce, as before, a system of included $i$-endings: $\ell_1 \subset \ell_2 \subset \ldots \subset \ell_n$ (we permute the letters of the alphabet if it is necessary). We show that the passage from $\ell_{i-1}$ to $\ell_i$ is possible by adding only two symbols (letters of alphabet {$\bf A$}). 

When we passed from $\ell_{i-1}$ to $\ell_i$ let there appear symbols $y$ and $z$. $\ell_{i-1}$ may be denoted by $AB$, where $A$ is a certain word, $B$ consists of the letters of the word $A$ (which are somehow mixed) and $B$ contains one letter $a_{i-1}$ less than $A$ does. Let $x$ be the last letter of the word $A$ on the right. Then $\ell_i$ may be denoted by $yzKxB$, where $A=Kx$.  From the definition of $\ell_i$ we have the equation
$${y} \cup {z} \cup K = B \cup {x} \cup {a_i}.$$
which from the definition of $K$ and $B$ is equivalent to 
$$2{x} \cup {a_i} = {y} \cup {z} \cup {a_i}.$$
It follows necessarily that $x=a_{i-1}$ and either $y=a_{i-1}$, $z=a_i$ or $y=a_i,z=a_{i-1}$. Suppose $y=a_{i-1}$, $z=a_i$. 

For example, we have the following crucial word for a 6-letter alphabet:
$$a_4a_5a_3a_4a_2a_3a_1a_2|a_6a_4a_3a_2a_1a_2a_3a_4a_6,$$
(the vertical line was drawn for a more convenient visual perception  of the word). 

This word is crucial and its length is equal to 17.

We consider a case of an arbitrary $n \ge 3$ defining the word $W$ as
$$W = a_{n-2}a_{n-1}a_{n-3}a_{n-2} \ldots a_1a_2|a_na_{n-2}a_{n-3} \ldots a_2a_1a_2 \ldots a_{n-3}a_{n-2}a_n.$$

Then $|W| = 2(n-2)+n-1+n-2=4n-7$. 

Let us verify that the word $W$ is crucial.  

If we write the letters $a_1$, $a_2$, $a_n$ to the right of the word $W$ we will obviously have prohibited subwords. Let $2<i<n$. Then if we write the letters $a_i$ we will have the prohibition $$a_{i-1}a_i \ldots a_1a_2a_na_{n-2} \ldots a_i|a_{i-1} \ldots a_2a_1a_2 \ldots a_{n-2}a_na_i,$$

since the composition vectors of the left and right subwords with respect to the vertical line are equal. 

Before proving that $W \in \widehat{{\bf S}_2^n}$ we make the following remark.

In the word $W$ we have $\ell_n \subset \ell_1 \subset \ldots \subset \ell_{n-2} \subset \ell_{n-1}$. Substituting $a_1$ for $a_n$, $a_2$ for $a_1$, \ldots , $a_{n}$ for $a_{n-1}$ we obtain another word 
$$U=a_{n-1}a_n \ldots a_2a_3|a_1a_{n-1} \ldots a_3a_2a_3 \ldots a_{n-1}a_1,$$
for which $\ell_1 \subset \ell_2 \subset \ldots \subset \ell_n$.

In both cases (before and after substitution of letters of the alphabet) we have the construction of a crucial word (which will be proved below) hence the same upper bound of the length of a minimal crucial word. 

For $W$ it is more convenient to show further that $W \in \widehat{{\bf S}_2^n}$.

We rewrite $W$ making in it the marks (1),(2), \ldots ,(2n-4), which number the gaps between letters of a word like this: $$(2n-4)a_{n-2}(2n-5)a_{n-1} \ldots (2)a_1(1)a_2|a_na_{n-2} \ldots a_2a_1a_2 \ldots a_{n-2}a_n.$$

In a possible prohibition we mark the left and right bounds. Note that the length of a prohibition is an even number, and each letter must occur an even number of times in a prohibition. The left bound of the prohibition must lie to the right of the mark (2n-5), since the letter $a_{n-1}$ enters $W$ once; 

It must lie to the left of the mark (1), since to the right of the mark (1) there is one letter $a_1$.

Note that if $m$ is even then $(m)$ is not the left bound of the possible prohibition. Actually in this case two variants are possible:

1) the prohibition does not cover the left letter $a_n$.

2) the prohibition covers the left letter $a_n$.

In the second case we have not a prohibition, since if the prohibition begins from the even mark, then it can not cover the second $a_n$.

In the first case the right bound of the prohibition lies to the left of $a_n$, hence the letter $a_{ \frac{m}{2}+1}$ enters the prohibition only once.

Suppose the prohibition begins from the mark $(m)$ and $m$ is odd.

There are two possible cases.

1) The prohibition does not cover the left letter $a_n$ (this case is impossible since the letter $a_{ \lfloor \frac{m}{2} \rfloor}$ occurs the prohibition once).

2) The prohibition covers the left $a_n$. Then it covers the right $a_n$ too, and the letter $a_{ \lfloor \frac{m}{2} \rfloor}$ occurs an odd number of times in the prohibition. So $W \in \widehat{{\bf S}_2^n}$  and hence $L_{min}({\bf S}_2^n) \le 4n-7$ for $n > 2$.

We give now a lower bound.

Since the length of a minimal crucial word must be odd, and the passage from $\ell_i$ to $\ell_{i+1}$ requires at least two letters, we have that a trivial lower bound of the length of a minimal crucial word is $2n-1$.

Let us now improve the lower bound. Obviously a minimal crucial word in which $\ell_1 \subset \ell_2 \subset \ldots \subset \ell_n$ has an even number of occurrences of the letter $a_i$ for $i=1, \ldots ,n-1$ and an odd number of occurrences of the letter $a_n$. The word $U$ has two letters $a_1$, two letters $a_2$, one letter $a_n$ and four of any other letter. From proposition~\ref{prop} we know that there does not exist a crucial word that has the fewer number of letters, hence the word $U$ gives us the lower bound of the length of a minimal crucial word. \end{proof} 

\section{The Set of Prohibitions ${\bf S}_3^{n, k}$}

\begin{theorem} We have $$L_{min}({\bf S}_3^{n, k})=2k+1.$$ \end{theorem}

\begin{proof} For the set of prohibitions ${\bf S}_3^{n, k}$ we must have $|A|=|B| \ge k+1$, where $AB$ is an arbitrary prohibition. So we have  

$$L_{min}({\bf S}_3^{n, k}) \ge 2k+1.$$

An upper bound is given by the construction $p_1p_2 \ldots p_kxp_1p_2 \ldots p_k$, where $x,p_i \in {\bf A}$, $i=1, \ldots ,k$ and $x \ne p_i$. \end{proof}

\begin{remark} The crucial word with respect to ${\bf S}_3^{1, k}$ is unique and its length is $2k+1$. \end{remark}

\begin{theorem} We have $$L_{max}({\bf S}_3^{2, k})=3k+3.$$ \end{theorem}

\begin{proof} Let 

\[ \bar{a} = \left\{ \begin{array}{ll}
1, & \mbox{if $a=2$,} \\
2, & \mbox{if $a=1$.}
\end{array}
\right. \] 

Moreover, let us consider an arbitrary crucial word $A$, with respect to ${\bf S}_3^{2, k}$, of length greater then $3k+3$. It is easy to see that if $a_1a_2 \ldots a_{k+1}$ are the first $k+1$ letters of $A$ then the next $k+1$ letters of $A$ must be $\bar{a}_1\bar{a}_2 \ldots \bar{a}_{k+1}$, because otherwise the first $2k+2$ letters of $A$ will form a prohibited subword. By the same argument, we can show that 
$$A=a_1a_2 \ldots a_{k+1}\bar{a}_1\bar{a}_2 \ldots \bar{a}_{k+1}a_1a_2 \ldots a_{k+1}\bar{a}_1 \ldots. $$ 

Let us consider the subwords $A_i$ of $A$ of the length $2k+4$ which start from the $i$th letter, where $1 \leq i \leq k$:
$$A_i=\underbrace{a_ia_{i+1} \ldots a_{k+1}\bar{a}_1 \ldots \bar{a}_i}_{k+2}\underbrace{\bar{a}_{i+1} \ldots \bar{a}_{k+1}a_1 \ldots a_{i+1}}_{k+2}$$ 
If $a_i=\bar{a}_{i+1}$ then the underbraced subwords of $A_i$ are the same in the first and in the last positions, so they differ in at most $k$ positions, hence $A_i$ is prohibited. So we must have $a_i=a_{i+1}$ for $i=1, \ldots, k$. 

Without loss of generality we can assume that $a_1=1$, so

$$A=\underbrace{11 \ldots 1}_{k+1}\underbrace{22 \ldots 2}_{k+1}\underbrace{11 \ldots 1}_{k+1}2 \ldots. $$ 

It is easy to see that if the length of $A$ is greater then $3k+3$ then $A$ has a prohibited subword of length $2k+4$:

$$A=\underbrace{11 \ldots 1}_{k}\overbrace{1\underbrace{22 \ldots 2}_{k+1}}\overbrace{\underbrace{11 \ldots 1}_{k+1}2}2 \ldots. $$ 
   
(here and then two braces above an word show us a disposition of a prohibited subword and, in particular, a disposition of parts of this subword that correspond to $X$ and $Y$ from the definition of the set of prohibitions ${\bf S}_3^{n, k}$).

So $L_{max}({\bf S}_3^{2, k}) \leq 3k+3$.

To prove the theorem it is sufficient to check that there are no prohibited subwords in the word $A=\underbrace{11 \ldots 1}_{k+1}\underbrace{22 \ldots 2}_{k+1}\underbrace{11 \ldots 1}_{k+1}$.

Obviously the left end of a possible prohibition can be only in the left block $\underbrace{ 1 \ldots 1 }_{k+1}$:

$$\overbrace{\underbrace{ 1 \ldots 1 }_{j} \underbrace{ 2 \ldots 2 }_{i}}^{} \overbrace{\underbrace{ 2 \ldots 2 }_{k-i+1} \underbrace{ 1 \ldots 1 }_{2i+j-k-1}}^{}$$

with 

\begin{equation} j+i \ge k+1 \end{equation}

Two cases are possible:

$$1.\ j \ge k-i+1$$
$$2.\ j < k-i+1$$

In the first case there is non-coincidence between the left and right parts of the prohibition in the first $k-i+1$ letters and in the last $i$ letters that is non-coincidence in $k+1$ letters. So this case is impossible.

In the second case we have non-coincidence in the first $j$ letters and in the last $2i+j-k-1$ letters. Hence we have non-coincidence in $2(i+j)-k-1$ letters, that according to (1) is greater than or equal to $k+1$.

It follows that the word $\underbrace{ 1 \ldots 1}_{k+1} \underbrace{ 2 \ldots 2 }_{k+1} \underbrace{ 1 \ldots 1 }_{k+1}$ does not contain a prohibition and thus the theorem is proved. \end{proof}

\begin{theorem}[Incompleteness] The set of prohibitions ${\bf S}_3^{n, k}$ for $n \ge 3$ is incomplete. \end{theorem} 

\begin{proof} Since the alphabet ${\bf A}$ is finite, there is no trivial solution of the problem (such as taking all letters of ${\bf A}$ and obtaining an infinite sequence with the properties needed). So to prove the incompleteness of the set ${\bf S}_3^{n, k}$ we have to show the existence of an infinite word which is free from the set of prohibitions ${\bf S}_3^{n, k}$. 

We consider the case $n=3$ and the alphabet ${\bf A}=\{ 1, 2, 3 \}$, since the incompleteness of the set of prohibitions ${\bf S}_3^{n, k}$ for the case $n > 3$ will follow from the incompleteness of the set of prohibitions for the case $n = 3$. 

Let ${\bf B}=\{ a, b, c \}$ be an alphabet.  ${\bf B}^{*}$ is the set of all words of the alphabet~${\bf B}$. 

We define the mapping $f$ as follows:  

$$\underbrace{ 1 \ldots 1 }_{k+1} \rightarrow a,\ \underbrace{ 2 \ldots 2 }_{k+1} \rightarrow b,\ \underbrace{3 \ldots 3}_{k+1} \rightarrow c.$$

The domain of the mapping $f$ is the set of words of the alphabet \ $${\bf C}=\{\ \underbrace{ 1 \ldots 1 }_{k+1},\ \underbrace{ 2 \ldots 2 }_{k+1},\ \underbrace{ 3 \ldots 3 }_{k+1}\ \}.$$ 

The image of the mapping $f$ is the set ${\bf B}^{*}$.

Let the set of prohibitions ${\bf S}^{\prime}=\{ XX | X \in {\bf B}^{*}\}$. Obviously, the set ${\bf S}^{\prime}$ coincides with the set ${\bf S}_2^n$ whenever ${\bf A} = {\bf B}$.

It is known \cite{1} that for the alphabet ${\bf B}$ there exists the infinite sequence $L^{\prime}$ which is free from the set of prohibitions ${\bf S}^{\prime}$. $L^{\prime}$ is built by iteration of morphisms:

$$a \rightarrow abc$$
$$b \rightarrow ac$$
$$c \rightarrow b$$

The morphism iteration procedure is as follows.

We start from the letter $a$. Then we substitute this letter with $abc$. Then we substitute each letter in $abc$ by the rule above. We obtain after this step $abcacb$. And so on. Executing this procedure an infinite number of times gives us the sequence $L^{\prime}$. 

Let us prove that the sequence $L=f^{-1}(L^{\prime})$ does not contain words prohibited by ${\bf S}_3^{3, k}$.

We are going to prove the statement by considering $L$ and all possible dispositions of words prohibited by ${\bf S}_3^{3, k}$.

The sequence $L$ is built up from the letters of the alphabet {\bf C} or in other words from the {\em blocks} $\underbrace{x \ldots x}_{k+1}$, where $x \in \{ 1, 2, 3 \}$. It means that there are only three different cases for a disposition of a possible prohibition in $L$.

{\bf Case 1.} $\overbrace{\underbrace{x \ldots x}_{k+1} \ldots \underbrace{y \ldots y}_{k+1}}\overbrace{\underbrace{z \ldots z}_{k+1} \ldots \underbrace{t \ldots t}_{k+1}}$\ ;

{\bf Case 2.} $\underbrace{x \ldots x}_{i}\overbrace{\underbrace{x \ldots x}_{k-i+1} \ldots \underbrace{y \ldots y}_{k+1}}\overbrace{\underbrace{z \ldots z}_{k+1} \ldots \underbrace{t \ldots t}_{k-i+1}}\underbrace{t \ldots t}_{i}$\ , where $0 < i < k+1$; 

{\bf Case 3.} $\underbrace{x \ldots x}_{i}\overbrace{\underbrace{x \ldots x}_{k-i+1} \ldots \underbrace{y \ldots y}_{\ell}}\overbrace{\underbrace{y \ldots y}_{k-\ell+1} \ldots \underbrace{t \ldots t}_{k-j+1}}\underbrace{t \ldots t}_{j}$\ , where $0 \le i,j,l \le k+1$.

Now we will consider these cases and show that each of them is impossible.

{\bf Case 1.} Let {\bf P} denote the prohibited subword (prohibition) under consideration, {\bf R} and {\bf L} denote the right and the left parts of {\bf P} respectively.

It is obvious that {\bf L} and  {\bf R} have the same number of blocks. Moreover, the $i$th block of {\bf L} (from the left to the right) is equal to the $i$th block of {\bf R}, because otherwise we have non-coincidence of {\bf L} and {\bf R} in at least $k+1$ letters which contradicts the fact that ${\bf P} \in {\bf S}_3^{3, k}$. So we have that ${\bf P}=WW$ for some $W \in {\bf C}^*$. 

Now, $f({\bf P})=f(W)f(W)$ is a subword of $L^{\prime}$. But $f(W)f(W) \in {\bf S}^{\prime}$ which is impossible by the properties of $L^{\prime}$. So Case 1 is impossible. 

We note that an important consequence of Case 1 is the following. If $\underbrace{x \ldots x}_{k+1}\underbrace{y \ldots y}_{k+1}$ is a subword of $L$ then $x \neq y$. 

{\bf Case 2.} If there are no letters between $\underbrace{x \ldots x}_{k+1}$ and $\underbrace{y \ldots y}_{k+1}$, that is
$${\bf P}=\overbrace{\underbrace{x \ldots x}_{k-i+1}\underbrace{y \ldots y}_{k+1}}\overbrace{\underbrace{z \ldots z}_{k+1}\underbrace{t \ldots t}_{k-i+1}}\ ,$$ 
then we must have $x=z$, because otherwise we have $x \neq z$ and $y \neq z$ which gives us that {\bf L} and {\bf R} differ in the first $k+1$ positions, but this contradicts ${\bf P} \in {\bf S}_3^{3, k}$.

By the same argument we have $y=t$, so 
$${\bf P}=\overbrace{\underbrace{x \ldots x}_{k-i+1}\underbrace{y \ldots y}_{k+1}}\overbrace{\underbrace{x \ldots x}_{k+1}\underbrace{y \ldots y}_{k-i+1}}\ .$$
But if we consider now $f(L)=L^{\prime}$ then it has 
$${\bf P}^{\prime}=\overbrace{f(\underbrace{x \ldots x}_{k+1})f(\underbrace{y \ldots y}_{k+1})}\overbrace{f(\underbrace{x \ldots x}_{k+1})f(\underbrace{y \ldots y}_{k+1})}\ .$$
as a subword, which is impossible since ${\bf P}^{\prime} \in {\bf S}^{\prime}$.

So there is some non-empty subword in {\bf L} between $\underbrace{x \ldots x}_{k+1}$ and $\underbrace{y \ldots y}_{k+1}$, and {\bf P} can be written as 
$${\bf P}=\overbrace{\underbrace{x \ldots x}_{k-i+1}\underbrace{x_1 \ldots x_1}_{k+1} \ldots \underbrace{x_p \ldots x_p}_{k+1}\underbrace{y \ldots y}_{k+1}}\overbrace{\underbrace{z \ldots z}_{k+1}\underbrace{z_1 \ldots z_1}_{k+1} \ldots \underbrace{z_p \ldots z_p}_{k+1}\underbrace{t \ldots t}_{k-i+1}}\ .$$
There are two possible subcases here.

{\bf 1.} $x=z$. 

Since $x \neq x_1$ we have $x_1 \neq z$. If $x_1 \neq z_1$ then {\bf L} and {\bf R} differ in $k+1$ position starting from the ($k-i+2$)th position, which is impossible since ${\bf P} \in {\bf S}_3^{3, k}$. So $x_1=z_1$.

In the same way, for each of $x_2$, $x_3$, \ldots $x_p$, $y$, we can obtain that
$${\bf P}=\overbrace{\underbrace{z \ldots z}_{k-i+1}\underbrace{z_1 \ldots z_1}_{k+1} \ldots \underbrace{z_p \ldots z_p}_{k+1}\underbrace{t \ldots t}_{k+1}}\overbrace{\underbrace{z \ldots z}_{k+1}\underbrace{z_1 \ldots z_1}_{k+1} \ldots \underbrace{z_p \ldots z_p}_{k+1}\underbrace{t \ldots t}_{k-i+1}}$$

which leads us to the fact that {\bf L} has a subword $WW$ for some $W \in {\bf C}^*$, hence ${\bf L}^{\prime}$ has a subword $f(W)f(W)$ which is impossible. 

So the subcase 1 is impossible.

{\bf 2.} $x \neq z$.

If $x_1 \neq z$ then {\bf L} and {\bf R} differ in $k+1$ position starting from the first position, which is impossible since ${\bf P} \in {\bf S}_3^{3, k}$. So $x_1=z$.

If $x_2 \neq z_1$ then {\bf L} and {\bf R} differ in $k+1$ position starting from the ($k+2$)th position, what is impossible by the same arguments as above. So $x_2=z_1$. And so on.

We have
$${\bf P}=\overbrace{\underbrace{x \ldots x}_{k-i+1}\underbrace{z \ldots z}_{k+1}\underbrace{z_1 \ldots z_1}_{k+1} \ldots \underbrace{z_p \ldots z_p}_{k+1}}\overbrace{\underbrace{z \ldots z}_{k+1}\underbrace{z_1 \ldots z_1}_{k+1} \ldots \underbrace{z_p \ldots z_p}_{k+1}\underbrace{t \ldots t}_{k-i+1}}\ .$$

Applying $f$ to $L$ gives us a subword ${\bf P}^{\prime}$ of $L^{\prime}$,
$${\bf P}^{\prime}=\overbrace{f(\underbrace{z \ldots z}_{k+1})f(\underbrace{z_1 \ldots z_1}_{k+1}) \ldots f(\underbrace{z_p \ldots z_p}_{k+1})}\overbrace{f(\underbrace{z \ldots z}_{k+1})f(\underbrace{z_1 \ldots z_1}_{k+1}) \ldots f(\underbrace{z_p \ldots z_p}_{k+1})},$$
which is prohibited in $L^{\prime}$ by ${\bf S}^{\prime}$. 

We have got that subcase 2 is impossible and hence Case 2 is impossible. 

{\bf Case 3.} We can assume that $\ell \neq 0$ and $\ell \neq k+1$, because otherwise we deal with either Case 1 or Case 2 which are impossible.

We suppose that $i \ge \ell$ (the case $i < \ell$ can be considered in the same way).

If there are no letters between $\underbrace{y \ldots y}_{k-\ell+1}$ and $\underbrace{t \ldots t}_{k-j+1}$, then we have either
$${\bf P}=\overbrace{\underbrace{x \ldots x}_{k-i+1}\underbrace{y \ldots y}_{\ell}}\overbrace{\underbrace{y \ldots y}_{k-\ell+1}\underbrace{t \ldots t}_{k-j+1}}$$ 
or
$${\bf P}=\overbrace{\underbrace{x \ldots x}_{k-i+1}\underbrace{z \ldots z}_{k+1}\underbrace{y \ldots y}_{\ell}}\overbrace{\underbrace{y \ldots y}_{k-\ell+1}\underbrace{t \ldots t}_{k-j+1}}\ .$$ 

In the first of these cases we have that $x \neq y$ and $y \neq t$ which gives us that {\bf L} and {\bf R} have non-coincidence in at least $k+1$ letters, but this contradicts ${\bf P} \in {\bf S}_3^{3, k}$.

In the second case we must have $z=t$, because otherwise since $z \neq y$ and $t \neq y$, {\bf L} and {\bf R} have non-coincidence in the last $k+1$ letters which is impossible. So in the second case we have
$${\bf P}=\overbrace{\underbrace{x \ldots x}_{k-i+1}\underbrace{t \ldots t}_{k+1}\underbrace{y \ldots y}_{\ell}}\overbrace{\underbrace{y \ldots y}_{k-\ell+1}\underbrace{t \ldots t}_{k-j+1}}\ .$$  

If $x \neq y$ then {\bf L} and {\bf R} have non-coincidence in the first $k-\ell+1$ positions and in the last $\ell$ positions, that is they have non-coincidence in at least $k+1$ positions which is impossible. So $x=y$.

Now applying $f$ to $L$ gives us that $L^{\prime}$ has a subword
$${\bf P}^{\prime}=\overbrace{f(\underbrace{x \ldots x}_{k+1})f(\underbrace{t \ldots t}_{k+1})}\overbrace{f(\underbrace{x \ldots x}_{k+1})f(\underbrace{t \ldots t}_{k+1})}$$

which is impossible.

So there is some non-empty subword in {\bf R} between $\underbrace{y \ldots y}_{k-\ell+1}$ and $\underbrace{t \ldots t}_{k-j+1}$, and {\bf P} can be written in the form 
$${\bf P}=\overbrace{\underbrace{x \ldots x}_{k-i+1}L_1 \ldots L_p\underbrace{y \ldots y}_{\ell}}\overbrace{\underbrace{y \ldots y}_{k-\ell+1}R_1 \ldots R_{p^{\prime}}\underbrace{t \ldots t}_{k-j+1}}\ ,$$
where $L_s,\ R_m \in {\bf C}$, for $1 \leq s \leq p$, $1 \leq m \leq p^{\prime}$, and either $p=p^{\prime}$ or $p=p^{\prime}+1$.   

We define $\bigtriangleup(L_s)=x_s$ if $L_s=\underbrace{x_s \ldots x_s}_{k+1}$. In the same way we define $\bigtriangleup(R_m)$.

Now we have that either $p=p^{\prime}$ or $p=p^{\prime} + 1$. Each of these cases has two possible subcases: either $x=y$ or $x \neq y$. Let us consider the case $p=p^{\prime} + 1$. The other case can be considered by similar reasoning. Thus we must consider the following subcases a) and b):

\ \ {\bf a)} $x=y$; It must be that $L_1=R_1$, because otherwise {\bf L} and {\bf R} differ in $k+1$ positions starting from the ($k-i+2$)th position. Then we consider one by one $L_2$, $L_3$, \ldots ,$L_p$. One can see that in this subcase
$${\bf P}=\overbrace{\underbrace{y \ldots y}_{k-i+1}R_1 \ldots R_{p^{\prime}}\underbrace{t \ldots t}_{k+1}\underbrace{y \ldots y}_{\ell}}\overbrace{\underbrace{y \ldots y}_{k-\ell+1}R_1 \ldots R_{p^{\prime}}\underbrace{t \ldots t}_{k-j+1}}\ ,$$ 
and $L$ has $WW$ as a subword, where $W=\underbrace{y \ldots y}_{k+1}R_1 \ldots R_{p^{\prime}}\underbrace{t \ldots t}_{k+1}$ which is impossible.

\ \ {\bf b)} $x \neq y$; There are two special subcases here, namely either $\bigtriangleup(L_1)=y$ or $L_1=R_1$.

If $\bigtriangleup(L_1)=y$ then 
$${\bf P}=\overbrace{\underbrace{x \ldots x}_{k-i+1}\underbrace{y \ldots y}_{k+1}R_1 \ldots R_{p^{\prime}}\underbrace{y \ldots y}_{\ell}}\overbrace{\underbrace{y \ldots y}_{k-\ell+1}R_1 \ldots R_{p^{\prime}}\underbrace{t \ldots t}_{k-j+1}}\ ,$$ 
and $L$ has $WW$ as a subword, where $W=\underbrace{y \ldots y}_{k+1}R_1 \ldots R_{p^{\prime}}$ which is impossible.

So $L_1=R_1$. In this case we have
$${\bf P}=\overbrace{\underbrace{x \ldots x}_{k-i+1}R_1 \ldots R_{p^{\prime}}\underbrace{t \ldots t}_{k+1}\underbrace{y \ldots y}_{\ell}}\overbrace{\underbrace{y \ldots y}_{k-\ell+1}R_1 \ldots R_{p^{\prime}}\underbrace{t \ldots t}_{k-j+1}}\ .$$ 

Since $y \neq x$, $y \neq \bigtriangleup(R_1)$ and $y \neq t$, {\bf L} and {\bf R} have non-coincidence in the first $k-l+1$ positions and in the last $l$ positions, so they have  non-coincidence in $k+1$ positions which contradicts ${\bf P} \in {\bf S}_3^{3, k}$.

We have got that Case 3 is impossible.

We have proved that the infinite word $L$ contains no word from the set ${\bf S}_3^{3, k}$ as a subword, therefore ${\bf S}_3^{n, k}$ is incomplete for $n \geq 3$. \end{proof}

\section{The Complexity of Problems on Completeness of Sets of Words}

It is known \cite{3,4} that the complexity of deciding whether or not an arbitrary set of prohibited words ${\bf S}$ is complete (or blocking) is $O(|{\bf S}| \cdot n)$, where $n$ is the greatest length of a word in ${\bf S}$. 

It is interesting in its own right to be able to effectively (in polynomial time) recognise whether a set is complete, but also to give a more detailed characterisation of the 
set of words $\widehat{{\bf S}}$, in particular to find the greatest length of a word that is free from ${\bf S}$. The set ${\bf A}^n$ is the set of all the words in the alphabet ${\bf A}$ whose length is equal to $n$. If ${\bf S} \subseteq {\bf A}^n$ and $L(n) = \max\limits_{\bf S} L(\widehat{{\bf S}})$, where $L(\widehat{{\bf S}})$ is the greatest length of a word that is free from ${\bf S}$, then \cite{3} we have
$$L(n) = {|{\bf A}|}^{n-1} + n - 2 = C(n) + n - 1.$$
Here $C(n)$ is the greatest length of a single path in the de Bruijn graph of order $n$ that has no chords and does not go through the vertices with loops corresponding to the constant words $(x, \ldots , x)$ where $x \in {\bf A}$.

One can find all words that are free from ${\bf S}$, in particular all crucial words, simply by considering all words of length less than or equal to $L(\widehat{{\bf S}})$ and checking for each word, if it is free from ${\bf S}$. Such an algorithm is not effective since it can require considering $|{\bf A}|^{L(n)}$ words.

The question of deciding the possible lengths of words that are free from ${\bf S}$, in particular of crucial words, can be formulated as a problem of recognising  properties of ``languages of prohibitions'' in the terminology of the theory of NP-completeness~\cite{6}. 

{\bf Problem A:} 

{\em Given:} An arbitrary set of words ${\bf S}$ and a natural number $\ell$.

{\em The question:} Does there exist a word of length at least $\ell$ that is free from ${\bf S}$?

In order to compare, we formulate the problem of completeness of a set of words ${\bf S}$ in the same form.

{\bf Problem B:} 

{\em Given:} An arbitrary set of words ${\bf S}$ in an alphabet ${\bf A}$.

{\em The question:} Does there exist $\ell \in {\bf N}$ such that $|X| \leq \ell$ for any word $X$ that is free from ${\bf S}$?

Considering problems A and B as problems of recognising properties of finite sets ${\bf S}$, we observe that problem B is a question of {\em existence of a bound} on the length of the words that are free from ${\bf S}$. This problem, as we have already mentioned, can be solved effectively with complexity of order $|{\bf S}| \cdot n$. In the same time the problem A is a question of {\em determining of this bound}. We will show that problem A, as opposed to problem B, is NP-complete.

The research on the problems of completeness of sets of words and languages of prohibited subwords was begun by different authors \cite{1,3,4,5,7,9} in the 1970s. The interest in the general question in this area arose from considerations of different types of special problems, in particular, in coding theory, combinatorics of symbolic sequences, number theory and problems of Ramsey type (for instance the arithmetic progressions in partitions of the natural row). For algebraic problems it is more typical to study avoidance of infinite sets ${\bf S}$ that are defined by prohibitions of words (called {\em terms}) in an alphabet of variables that can themselves be words \cite{1,9}. Different problems on sequences without repetitions, under variation the concept of ``strong'' or ``weak'' repetition of subwords, are the typical examples of problems of this class. Finally we observe that problems A and B for infinite sets ${\bf S}$ do not make sense if one does not consider particular constructive methods for generating a set~${\bf S}$. 

Let ${\bf A}=\{ a_1, \ldots ,a_n\}$ be an alphabet and ${\bf A}_{\ell}$ be the set of all those words on the alphabet ${\bf A}$ whose length is less than or equal to $\ell$. We assume also that the empty word belongs to ${\bf A}_{\ell}$ and that ${\bf S}_1$ is an arbitrary set such that ${\bf S}_1 \subseteq {\bf A}_2 \setminus {\bf A}_1$. We define ${\bf S}_2$ by 
$${\bf S}_2 = \{ xXx |\ x \in {\bf A}, X \in {\bf A}^{n-1} \}.$$

So the set ${\bf S}_2$ contains all possible words of length less than or equal to $n+1$ whose first letter coincides with their last letter. Suppose ${\bf S} = {\bf S}_1 \cup {\bf S}_2$.

We now consider an ``auxiliary'' problem $\mbox{A}^{\prime}$. 

{\bf Problem $\mbox{A}^{\prime}$:}

{\em Given:} A set ${\bf S}$ of the type described above and a natural number $\ell$, $\ell \leq n$.

{\em The question:} Does there exist a word of length at least $\ell$ that is free from ${\bf S}$?

In case of the problem $\mbox{A}^{\prime}$, the restriction on $\ell$ is natural, because any word free from ${\bf S}$ is free from ${\bf S}_2$ and therefore consists of different letters of the alphabet, whence its length is less than or equal to $n$.

Checking whether a given word of length $\ell$ (a solution of $\mbox{A}^{\prime}$ that we ``guessed'') is free from ${\bf S}$ can be done in polynomial time. Indeed, the freeness from ${\bf S}_2$ of the word is equivalent to the absence of identical letters in the word (which can be checked in linear time) and the freeness from ${\bf S}_1$ is recognised by considering all subwords of length 2 (there are $\ell -1$ such subwords) and by checking for each of them whether it belongs to ${\bf S}_1$ (polynomial checking time).

We now introduce the problem of ``the longest path in a graph'', which is known to be NP-complete (see \cite{6}).

{\bf Problem ``path'':}

{\em Given:} A directed graph $\vec{G}(V,E)$ and a natural number $\ell$, $\ell \leq |V|=n$.

{\em The question:} Does there exist a simple directed path (without self-intersections in vertices) of length at least $\ell$? 

One can obtain a correspondence between problem $\mbox{A}^{\prime}$ and problem ``path'' as follows. We compare vertices $v_1, \ldots, v_n$ from $V(\vec{G})$ to the letters $a_1, \ldots, a_n$ in the alphabet ${\bf A}$. Also we compare each edge $\vec{v_iv_j}$ from $E(\vec{G})$ to the word $a_ia_j$. We form the set ${\bf S}_1$ from all such words of ${\bf A}^2$ that correspond to the edges of the graph that is the complement of $\vec{G}$ with respect to the complete directed graph.

Now to any oriented simple path $v_{i_1}, \ldots ,v_{i_{\ell}}$ of length $\ell$ in $\vec{G}$ there corresponds the word $a_{i_1} \ldots a_{i_{\ell}}$ of length $\ell$, consecutive letters of which correspond to vertices in the order in which the path passed through them. This word is free from ${\bf S}_1$ because $a_{i_j}a_{i_{j+1}} \not\in {\bf S}_1$ for any $i=1,2, \ldots ,\ell -1$. The word is free from the set ${\bf S}_2$ as well because in the path there is no repetition of vertices (a property of a simple path) and therefore $a_{i_1} \ldots a_{i_{\ell}}$ does not contain a subword of the form $a_iXa_i$ for any word $X$ and any letter $a_i \in {\bf A}$.    

Conversely, to any word in the alphabet ${\bf A}$ that is free from ${\bf S}$ there corresponds a path in $\vec{G}(V,E)$ that goes through edges from $E(\vec{G})$ since the word is free from ${\bf S}_1$ and that is not self-intersected since the word is free from ${\bf S}_2$.

Now NP-completeness of problem $\mbox{A}^{\prime}$ and the more general problem A follows from NP-completeness of the problem ``path''.


\begin{thebibliography}{9}
\bibitem{1} Bean D. R., Ehrenfeucht A., McNulty G. F., Avoidable patterns in strings of symbols, Pacific J. Math., Vol. {\bf 85}, No. 2, (1979), 261-294.
\bibitem{2} Choffrut C., Karhum\"aki J., {\em Handbook of formal languages}, Vol. {\bf 1}: Word, language, grammar, Berlin, Springer, (1997), 329-438. 
\bibitem{3} Evdokimov A., Complete sets of words and their numerical characteristics, Metody Diskret. Analiz., Novosibirsk, IM SB RAS, No. {\bf 39}, (1983), 7-32 (in Russian). \\
Also 86e:68087 in the Mathematical Reviews on the Web.
\bibitem{4} Evdokimov A., The completeness of sets of words, Proceedings of the All-Union seminar on discrete mathematics and its applications (Russian) (Moscow, 1984), Moskov. Gos. Univ., Mekh.-Mat. Fak., Moscow, (1986), 112-116 (in Russian). \\
Also 89e:68066 in the Mathematical Reviews on the Web.
\bibitem{5} Evdokimov A., Krainev V., Problems on completeness of sets of words, Proceedings of the 22nd regional scientific conference by the Popov association, Novosibirsk, (1979), 105-107 (in Russian).
\bibitem{6} Garey M., Johnson D., {\em Computers and intractability: a guide to the theory of NP-completeness}, W. H. Freeman, (1979).
\bibitem{7} Lothaire M., {\em Combinatorics on Words}, Encyclopedia of Mathematics, Vol. {\bf 17}, Addison-Wesley (1986). 
\bibitem{8} Salomaa A., {\em Jewels of Formal Language Theory}, Computer Science Press, (1981).
\bibitem{9} Zimin A. I., Blocking sets of terms, Mat. Sbornik, Vol. {\bf 119}, No. 3, 363-375, (1982) (in Russian). \\
Also 84d:20072 in the Mathematical Reviews on the Web.
\end{thebibliography}

\end{document}